\documentclass[11pt]{amsart}
\usepackage{amsmath,amssymb,mathrsfs}
\newtheorem{theorem}{Theorem}[section]
\newtheorem{corollary}[theorem]{Corollary}
\newtheorem{conjecture}[theorem]{Conjecture}
\newtheorem*{unif.thm}{Uniformization Theorem}
\newtheorem*{yamabe.conj}{Yamabe Conjecture}
\newtheorem*{diffeo.sph.thm}{Differentiable Sphere Theorem}
\theoremstyle{definition}
\newtheorem{definition}[theorem]{Definition}

\begin{document}

\title{Evolution equations in Riemannian geometry}
\author{Simon Brendle}
\address{Department of Mathematics \\ Stanford University \\ Stanford, CA 94305}
\begin{abstract}
A fundamental question in Riemannian geometry is to find canonical metrics on a given smooth manifold. In the 1980s, R.~Hamilton proposed an approach to this question based on parabolic partial differential equations. The goal is to start from a given initial metric and deform it to a canonical metric by means of an evolution equation. There are various natural evolution equations for Riemannian metrics, including the Ricci flow and the conformal Yamabe flow. In this survey, we discuss the global behavior of the solutions to these equations. In particular, we describe how these techniques can be used to prove the Differentiable Sphere Theorem.

This article is based on the Takagi Lectures delivered by the author at the Research Institute for Mathematical Sciences, Kyoto University, on June 4, 2011.
\end{abstract}

\maketitle 

\section{The Ricci flow on surfaces and the uniformization theorem} 

One of the central questions in modern differential geometry is concerned with the existence of canonical metrics on a given manifold $M$. This question has a long history (going back to ideas of H.~Hopf and R.~Thom), and is inspired by the uniformization theorem in dimension $2$. In fact, if $M$ is a two-dimensional surface, then the uniformization theorem implies that $M$ admits a metric of constant curvature:

\begin{unif.thm}
Let $M$ be a compact surface, equipped with a Riemannian metric $g_0$. Then there exists a Riemannian metric $g$ on $M$ which is conformal to $g_0$ and has constant scalar curvature.
\end{unif.thm}

Recall that two metrics $g_0$ and $g$ are said to be conformal if there exists a function $w \in C^\infty(M)$ such that $g = e^{2w} \, g_0$.

The Ricci flow provides a canonical deformation from the given metric $g_0$ to a constant curvature metric. To describe this, let $M$ be a compact surface, equipped with an initial metric $g_0$. We then evolve the metric by the evolution equation 
\begin{equation} 
\label{2d.ricci.flow}
\frac{\partial}{\partial t} g(t) = -R_{g(t)} \, g(t), 
\end{equation}
where $R_{g(t)}$ denotes the scalar curvature of the evolving metric $g(t)$. The evolution equation (\ref{2d.ricci.flow}) was first studied by R.~Hamilton \cite{Hamilton3}. Note that 
\begin{equation} 
\frac{\partial}{\partial t} R_{g(t)} = \Delta_{g(t)} R_{g(t)} + R_{g(t)}^2 
\end{equation}
for any solution $g(t)$ of (\ref{2d.ricci.flow}).

Since the evolution equation (\ref{2d.ricci.flow}) preserves the conformal structure, we may write the metric $g(t)$ in the form $g(t) = e^{2w(t)} \, g_0$ for some function $w(t) \in C^\infty(M)$. The scalar curvature of $g(t)$ is given by the formula 
\[R_{g(t)} = e^{-2w(t)} \, \Big ( -2 \, \Delta_{g_0} w(t) + R_{g_0} \Big ).\] 
Substituting this formula into (\ref{2d.ricci.flow}) leads to the following evolution equation for the conformal factor: 
\begin{equation} 
\label{conf.factor}
\frac{\partial}{\partial t} w(t) = e^{-2w(t)} \, \Big ( \Delta_{g_0} w(t) - \frac{1}{2} \, R_{g_0} \Big ). 
\end{equation} 
Note that the function $w(t)$ satisfies a nonlinear partial differential equation of parabolic type; in particular, the evolution equation (\ref{2d.ricci.flow}) always admits a shorttime solution. Moreover, the equation (\ref{conf.factor}) is closely related to the porous medium equation; see e.g. \cite{Daskalopoulos-Hamilton}.

The longtime behavior of the flow depends in a crucial way on the sign of the Euler characteristic of $M$. In fact, the Gauss-Bonnet theorem implies that 
\[\frac{d}{dt} \text{\rm vol}(M,g(t)) = -\int_M R_{g(t)} \, d\text{\rm vol}_{g(t)} = -4\pi\chi(M),\] 
where $\chi(M)$ denotes the Euler characteristic of $M$. In particular, if $\chi(M)$ is negative, then the area of $M$ will increase; on the other hand, the area of $M$ will decrease if $\chi(M)$ is positive.

In the seminal paper \cite{Hamilton3}, Hamilton analyzed the global behavior of the solutions to (\ref{2d.ricci.flow}).

\begin{theorem}[R.~Hamilton \cite{Hamilton3}]
\label{convergence.of.2d.ricci.flow}
Let $M$ be a compact surface, equipped with a Riemannian metric $g_0$. Moreover, let $g(t)$, $t \in [0,T)$, be the unique maximal solution to the Ricci flow with initial metric $g_0$. 
\begin{itemize}
\item[(i)] If $\chi(M) < 0$, then $T = \infty$, and the rescaled metrics $\frac{1}{t} \, g(t)$ converge to a metric of constant scalar curvature $-1$ as $t \to \infty$.
\item[(ii)] If $\chi(M) = 0$, then $T = \infty$, and the metrics $g(t)$ converge to a flat metric as $t \to \infty$.
\item[(iii)] If $g_0$ has positive scalar curvature, then 
\[T = \frac{\text{\rm vol}(M,g_0)}{4\pi\chi(M)}.\] 
Moreover, the rescaled metrics $\frac{1}{T-t} \, g(t)$ converge to a metric of constant scalar curvature $1$ as $t \to T$.
\end{itemize}
\end{theorem}

In statement (iii), the assumption that $g_0$ has positive scalar curvature can be replaced by the weaker assumption that the Euler characteristic $\chi(M)$ is positive (cf. \cite{Chow1}).

The case when $g_0$ has positive scalar curvature is the most difficult one. In order to handle this case, Hamilton proved a remarkable monotonicity formula. In fact, Hamilton showed that the function 
\[t \mapsto \int_M R_{g(t)} \, \log \big ( (T-t) \, R_{g(t)} \big ) \, d\text{\rm vol}_{g(t)}\] 
is monotone decreasing. Using this monotonicity formula and a blow-up argument, it follows that the product $(T-t) \, R_{g(t)}$ is uniformly bounded from above. Hence, there exists a sequence of times $t_k \to T$ with the property that the rescaled metrics $\frac{1}{T-t_k} \, g(t_k)$ converge to a Ricci soliton metric in the Cheeger-Gromov sense. On the other hand, Hamilton showed that any compact two-dimensional Ricci soliton has constant curvature (see also \cite{Chen-Lu-Tian} for a simplified argument, which avoids the use of the uniformization theorem). Consequently, the rescaled metrics $\frac{1}{T-t_k} \, g(t_k)$ must converge to a metric of constant scalar curvature.

M.~Struwe \cite{Struwe2} has given an alternative proof of Theorem \ref{convergence.of.2d.ricci.flow} using PDE techniques. This approach relies on a concentration-compactness argument, as developed by X.~Chen \cite{Chen2}, \cite{Chen3} in the context of the Calabi flow.

\section{The conformal Yamabe flow}

In 1960, H.~Yamabe \cite{Yamabe} proposed a generalization of the uniformization theorem to higher dimensions. More precisely, Yamabe conjectured the following: 

\begin{yamabe.conj}
Let $M$ be a compact manifold of dimension $n \geq 3$, and let $g_0$ be a Riemannian metric on $M$. Then there exists a metric $g$ on $M$ which is conformal to $g_0$ and has constant scalar curvature.
\end{yamabe.conj}

Yamabe's Conjecture was confirmed by the work of N.~Trudinger \cite{Trudinger}, T.~Aubin \cite{Aubin}, and R.~Schoen \cite{Schoen1}. Building upon earlier work of Trudinger, Aubin was able to prove Yamabe's Conjecture when $n \geq 6$ and $(M,g_0)$ is not locally conformally flat. The remaining cases were settled by Schoen using the positive mass theorem. A.~Bahri \cite{Bahri} later gave an alternative proof in the locally conformally flat case.

In the 1980s, R.~Hamilton proposed a heat flow approach to the Yamabe problem. To describe this, let $g(t)$ be a smooth one-parameter family of Riemannian metrics on $M$. We say that $g(t)$ is a solution of the unnormalized Yamabe flow if 
\begin{equation}
\label{unnormalized.yamabe.flow}
\frac{\partial}{\partial t} g(t) = -R_{g(t)} \, g(t),
\end{equation} 
where $R_{g(t)}$ denotes the scalar curvature of $g(t)$. 

It is often convenient to consider a normalized version of the flow. We say that $g(t)$ is a solution of the normalized Yamabe flow if 
\begin{equation}
\label{normalized.yamabe.flow}
\frac{\partial}{\partial t} g(t) = -(R_{g(t)} - r_{g(t)}) \, g(t).
\end{equation} 
Here, $r_{g(t)}$ denotes the mean value of the scalar curvature of the metric $g(t)$; that is, 
\[r_{g(t)} = \frac{\int_M R_{g(t)} \, d\text{\rm vol}_{g(t)}}{\text{\rm vol}(M,g(t))}.\] 
The evolution equations (\ref{unnormalized.yamabe.flow}) and (\ref{normalized.yamabe.flow}) are equivalent in the sense that any solution of the equation (\ref{unnormalized.yamabe.flow}) can be transformed into a solution of (\ref{normalized.yamabe.flow}) by a rescaling procedure. In the sequel, we will always consider the normalized Yamabe flow (\ref{normalized.yamabe.flow}). Note that 
\begin{equation} 
\frac{\partial}{\partial t} R_{g(t)} = (n-1) \, \Delta_{g(t)} R_{g(t)} + R_{g(t)} \, (R_{g(t)} - r_{g(t)}) 
\end{equation}
if $g(t)$ is a solution of (\ref{normalized.yamabe.flow}).

Since the Yamabe flow preserves the conformal structure, we may write $g(t) = u(t)^{\frac{4}{n-2}} \, g_0$ for some fixed background metric $g_0$. The scalar curvature of the metric $g(t)$ is given by the formula 
\[R_{g(t)} = u(t)^{-\frac{n+2}{n-2}} \, \Big ( -\frac{4(n-1)}{n-2} \, \Delta_{g_0} u(t) + R_{g_0} \, u(t) \Big ).\] 
Therefore, the evolution equation (\ref{normalized.yamabe.flow}) is equivalent to the following evolution equation for the conformal factor: 
\begin{equation} 
\label{pde.for.conf.factor}
\frac{\partial}{\partial t} (u(t)^{\frac{n+2}{n-2}}) = \frac{n+2}{4} \, \Big ( \frac{4(n-1)}{n-2} \, \Delta_{g_0} u(t) - R_{g_0} \, u(t) + r_{g(t)} \, u(t)^{\frac{n+2}{n-2}} \Big ). 
\end{equation} 
As above, the conformal factor satisfies a nonlinear partial differential equation of parabolic type. The Yamabe flow can be viewed as the gradient flow for the Yamabe functional $E_{g_0}$. This functional is defined by 
\[E_{g_0}(u) = \frac{\int_M \big ( \frac{4(n-1)}{n-2} \, |du|_{g_0}^2 + R_{g_0} \, u^2 \big ) \, d\text{\rm vol}_{g_0}}{\big ( \int_M u^{\frac{2n}{n-2}} \, d\text{\rm vol}_{g_0} \big )^{\frac{n-2}{n}}}.\] 
The functional $E_{g_0}(u)$ is also referred to as the Yamabe energy. 

We note that the Yamabe functional $E_{g_0}(u)$ is closely related to the Sobolev quotient associated with the embedding of $H^1(M,g_0)$ into $L^{\frac{2n}{n-2}}(M,g_0)$. To fix notation, we denote by 
\[Y(M,g_0) = \inf_{0 < u \in C^\infty(M)} E_{g_0}(u)\] 
the infimum of the Yamabe functional of $(M,g_0)$.

The Yamabe functional also has a natural geometric interpretation. In fact, for any positive function $u \in C^\infty(M)$ we have
\[E_{g_0}(u) = \frac{\int_M R_g \, d\text{\rm vol}_g}{\text{\rm vol}(M,g)^{\frac{n-2}{n}}},\] 
where $g = u^{\frac{4}{n-2}} \, g_0$. In other words, the Yamabe functional $E_{g_0}$ may be viewed as the restriction of the normalized Einstein-Hilbert action to the conformal class of $g_0$. 

It is an interesting question to study the longtime behavior of the Yamabe flow. It was shown by R.~Hamilton \cite{Hamilton-notes} that the Yamabe flow (\ref{normalized.yamabe.flow}) cannot develop singularities in finite time. More precisely, Hamilton proved the following theorem:

\begin{theorem}[R.~Hamilton]
\label{global.existence}
Let $M$ be a compact manifold of dimension $n \geq 3$. Given any initial metric $g_0$ on $M$, the Yamabe flow (\ref{normalized.yamabe.flow}) admits a solution which is defined for all $t \geq 0$.
\end{theorem}

Once we know that the Yamabe flow admits a global solution, one would like to understand the asymptotic behavior as $t \to \infty$. In this direction, Hamilton proposed the following conjecture: 

\begin{conjecture}[R.~Hamilton]
Let $(M,g_0)$ be a compact Riemannian manifold of dimension $n \geq 3$, and let $g(t)$ be the unique solution of the Yamabe flow with initial metric $g_0$. Then the metrics $g(t)$ converge to a metric of constant scalar curvature as $t \to \infty$.
\end{conjecture}

Hamilton's Conjecture was first taken up by B.~Chow \cite{Chow2}. Chow showed that if $(M,g_0)$ is locally conformally flat and has positive Ricci curvature, then the Yamabe flow evolves $g_0$ to a metric of constant scalar curvature. This result was subsequently improved by R.~Ye \cite{Ye}.

\begin{theorem}[R.~Ye \cite{Ye}]
\label{lcf.case}
Let $(M,g_0)$ be a compact Riemannian manifold of dimension $n \geq 3$, and let $g(t)$ be the unique solution of the Yamabe flow with initial metric $g_0$. If $(M,g_0)$ is locally conformally flat, then the metrics $g(t)$ converge to a metric of constant scalar curvature as $t \to \infty$.
\end{theorem}

The proof of Theorem \ref{lcf.case} uses the developing map in a crucial way. More precisely, since $(M,g_0)$ is locally conformally flat, there exists a conformal diffeomorphism from the universal cover of $(M,g_0)$ onto a dense domain $\Omega \subset S^n$ (cf. \cite{Schoen-Yau2}). In particular, the metric $g(t)$ induces a metric $\tilde{g}(t)$ on $\Omega$. The metrics $\tilde{g}(t)$ are conformal to the standard metric on $S^n$, and evolve by the Yamabe flow. By applying the method of moving planes, Ye was able to prove a gradient estimate for the conformal factor. As a consequence, Ye obtained uniform upper and lower bounds for the conformal factor. A similar argument was used by R.~Schoen to prove a-priori estimates for solutions of the elliptic Yamabe equation (see \cite{Schoen2}).

Schwetlick and Struwe \cite{Schwetlick-Struwe} were able to remove the assumption that $(M,g_0)$ is locally conformally flat; instead, they assumed that the Yamabe energy of the initial metric is below a certain critical threshold: 

\begin{theorem}[H.~Schwetlick, M.~Struwe \cite{Schwetlick-Struwe}]
\label{schwetlick.struwe.thm}
Let $(M,g_0)$ be a compact Riemannian manifold of dimension $n$, where $3 \leq n \leq 5$. Moreover, let $g(t)$ be the unique solution to the Yamabe flow with initial metric $g_0$. If the Yamabe energy of $g_0$ is less than $\big [ Y(M,g_0)^{\frac{n}{2}} + Y(S^n)^{\frac{n}{2}} \big ]^{\frac{2}{n}}$, then the metrics $g(t)$ converge to a metric of constant scalar curvature as $t \to \infty$.
\end{theorem}

Here, $Y(S^n)$ denotes the infimum of the Yamabe functional on the round sphere $S^n$.

In \cite{Brendle1}, we removed the assumption on the Yamabe energy of $g_0$, and proved the convergence of the Yamabe flow for arbitrary initial metrics (see also \cite{Brendle2}):

\begin{theorem}[S.~Brendle \cite{Brendle1}]
\label{convergence.yamabe.flow}
Let $(M,g_0)$ be a compact Riemannian manifold of dimension $n$. We assume that either $3 \leq n \leq 5$ or $(M,g_0)$ is locally conformally flat. Moreover, let $g(t)$, $t \geq 0$, be the unique solution to the Yamabe flow with initial metric $g_0$. Then the metrics $g(t)$ converge to a metric of constant scalar curvature as $t \to \infty$.
\end{theorem}

Let us sketch the main ideas involved in the proofs of Theorem \ref{schwetlick.struwe.thm} and Theorem \ref{convergence.yamabe.flow}. As above, we write $g(t) = u(t)^{\frac{4}{n-2}} \, g_0$, where $u(t) \in C^\infty(M)$ is a positive function on $M$. If the conformal factor $u(t)$ is uniformly bounded from above, then the convergence of the flow follows from a general theorem due to L.~Simon \cite{Simon}. Hence, it suffices to show that the function $u(t)$ is uniformly bounded. To that end, one performs a blow-up analysis. Suppose that $\sup_M u(t_k) \to \infty$ for some sequence of times $t_k \to \infty$. It follows from a result of Struwe \cite{Struwe1} that the flow forms a finite number of spherical bubbles. 

In the setting of Theorem \ref{schwetlick.struwe.thm}, the assumption on the initial energy implies that the flow forms at most one bubble. Moreover, Schwetlick and Struwe were able to rule out this possibility using the positive mass theorem of Schoen and Yau \cite{Schoen-Yau1}. The case of higher initial energy is much more subtle, as the flow might form multiple bubbles, and the interactions between different bubbles must be taken into account. To rule out this possibility, we again use the positive mass theorem; see \cite{Brendle1} for a detailed proof.

Finally, we describe a convergence result for the Yamabe flow in dimension $n \geq 6$. To fix notation, let $(M,g_0)$ be a compact Riemannian manifold of dimension $n \geq 6$. We denote by $\mathcal{Z}$ the set of all points $p \in M$ such that 
\[\limsup_{x \to p} d(p,x)^{-[\frac{n-6}{2}]} \, |W_{g_0}(x)| = 0,\] 
where $W_{g_0}$ denotes the Weyl tensor of the background metric $g_0$ and $d(\cdot,\cdot)$ denotes the Riemannian distance. Note that the set $\mathcal{Z}$ depends only on the conformal class of $g_0$.

\begin{theorem}[S.~Brendle \cite{Brendle3}]
\label{higher.dim}
Let $(M,g_0)$ be a compact Riemannian manifold of dimension $n \geq 6$. We assume that either $M$ is spin or $\mathcal{Z} = \emptyset$. Moreover, let $g(t)$, $t \geq 0$, be the unique solution to the Yamabe flow with initial metric $g_0$. Then the metrics $g(t)$ converge to a metric of constant scalar curvature as $t \to \infty$.
\end{theorem}

Note that the limit metric need not be a global minimizer of the Yamabe functional. In fact, there may be infinitely many metrics of constant scalar curvature in a given class (see e.g. \cite{Brendle4}, \cite{Brendle-Marques}).

\section{The Ricci flow in higher dimensions}

In this section, we discuss another natural evolution equation for Riemannian metrics. Let $M$ be a compact manifold of dimension $n \geq 3$, and let $g(t)$ be a smooth one-parameter family of Riemannian metrics on $M$. We say that the metrics $g(t)$ evolve by the Ricci flow if 
\begin{equation} 
\label{ricci.flow}
\frac{\partial}{\partial t} g(t) = -2 \, \text{\rm Ric}_{g(t)}, 
\end{equation}
where $\text{\rm Ric}_{g(t)}$ denotes the Ricci tensor of the metric $g(t)$. This evolution was introduced in 1982 in a landmark paper by R.~Hamilton \cite{Hamilton1}. Hamilton also considered a normalized version of the Ricci flow, which keeps the volume of $M$ fixed. In what follows, we will restrict ourselves to the unnormalized Ricci flow (\ref{ricci.flow}).

The evolution equations (\ref{unnormalized.yamabe.flow}) and (\ref{ricci.flow}) can both be viewed as generalizations of the Ricci flow in dimension $2$. Unlike the Yamabe flow, the Ricci flow does not preserve the conformal structure.

In \cite{Hamilton1}, Hamilton proved a shorttime existence theorem for the Ricci flow: 

\begin{theorem}[R.~Hamilton \cite{Hamilton1}]
\label{shorttime.existence}
Let $M$ be a compact manifold, equipped with a Riemannian metric $g_0$. Then there exists a real number $T > 0$ and a solution $g(t)$, $t \in [0,T)$, of the Ricci flow with the property that $g(0) = g_0$. Moreover, the solution $g(t)$ is unique.
\end{theorem}

The question of shorttime existence is subtle because the Ricci flow is only weakly parabolic. To overcome this obstacle, Hamilton employed the Nash-Moser inverse function theorem. A simplified proof was later given by D.~DeTurck \cite{DeTurck}.

Hamilton \cite{Hamilton1}, \cite{Hamilton2} studied how the curvature tensor changes when the metric is evolved by the Ricci flow. Hamilton showed that the Riemannian curvature tensor satisfies the following nonlinear heat equation: 
\begin{equation} 
\label{evolution.of.curvature} 
D_t R_{ijkl} = \Delta R_{ijkl} + \sum_{p,q} R_{ijpq} \, R_{klpq} + 2 \sum_{p,q} (R_{ipkq} \, R_{jplq} - R_{iplq} \, R_{jpkq}). 
\end{equation}
Here, $R$ denotes the Riemann curvature tensor of the metric $g(t)$, $\Delta$ denotes the Laplace operator of $g(t)$, and $D_t$ represents the covariant time derivative (cf. \cite{Brendle-book}, Section 2.3).

Using the evolution equation (\ref{evolution.of.curvature}) and the maximum principle, Hamilton \cite{Hamilton1} showed that the Ricci flow preserves positive Ricci curvature in dimension $3$. More precisely, let $(M,g_0)$ be a compact three-manifold, and let $g(t)$, $t \in [0,T)$, be the unique solution to the Ricci flow with initial metric $g_0$. If $(M,g_0)$ has positive Ricci curvature, then $(M,g(t))$ will have positive Ricci curvature for all $t \in [0,T)$. Moreover, Hamilton was able to give a complete description of the longtime behavior of the Ricci flow in this situation.

\begin{theorem}[R.~Hamilton \cite{Hamilton1}]
\label{Hamilton.dim.3}
Let $(M,g_0)$ be a compact three-manifold with positive Ricci curvature, and let $g(t)$, $t \in [0,T)$, denote the unique maximal solution to the Ricci flow with initial metric $g_0$. Then $T < \infty$, and the rescaled metrics $\frac{1}{4(T-t)} \, g(t)$ converge to a metric of constant sectional curvature $1$ as $t \to T$.
\end{theorem} 

In particular, if a three-manifold $M$ admits a metric of positive Ricci curvature, then $M$ admits a metric of constant sectional curvature. This confirmed an earlier conjecture of J.P.~Bourguignon \cite{Bourguignon}. In particular, Hamilton was able to draw the following topological conclusion:

\begin{corollary}[R.~Hamilton \cite{Hamilton1}]
\label{cor}
Let $M$ be a compact three-manifold which admits a metric of positive Ricci curvature. Then $M$ is diffeomorphic to a quotient $S^3 / \Gamma$, where $\Gamma$ is a finite group of standard isometries acting freely.
\end{corollary}

The quotient manifolds $S^3 / \Gamma$ are completely classified (cf. \cite{Wolf}). Hence, Corollary \ref{cor} provides a classification of all three-manifolds which admit metrics with positive Ricci curvature.

Theorem \ref{Hamilton.dim.3} has inspired a large body of work over the past 30 years. In particular, two lines of inquiry have been pursued:

First, it is very interesting to study the global properties of the Ricci flow in dimension $3$ for general initial metrics (i.e. without any curvature assumptions). This is a very subtle problem, as the Ricci flow may develop singularities. In \cite{Hamilton-survey}, Hamilton outlined a program for understanding these singularities. In particular, Hamilton proposed a surgery procedure in order to extend the flow past singularities (cf. \cite{Hamilton-survey}, \cite{Hamilton5}). This program culminated in G.~Perelman's proof of the Poincar\'e Conjecture (cf. \cite{Perelman1}, \cite{Perelman2}, \cite{Perelman3}). Non-technical surveys can be found in \cite{Besson}, \cite{Ecker}, or \cite{Leeb}. 

Second, one would like to generalize the convergence theory for the Ricci flow to dimensions greater than $3$. This requires some assumptions on the curvature of the initial metric. Recall that a Riemannian manifold $M$ is said to have positive curvature operator if 
\[\sum_{i,j,k,l} R_{ijkl} \, \varphi^{ij} \, \varphi^{kl} > 0\] 
for each point $p \in M$ and every non-zero two-form $\varphi \in \wedge^2 T_p M$. Moreover, $M$ is said to have two-positive curvature operator if 
\[\sum_{i,j,k,l} R_{ijkl} \, (\varphi^{ij} \, \varphi^{kl} + \psi^{ij} \, \psi^{kl}) > 0\] 
for all points $p \in M$ and all two-forms $\varphi,\psi \in \wedge^2 T_p M$ satisfying $|\varphi|^2 = |\psi|^2 = 1$ and $\langle \varphi,\psi \rangle = 0$.

R.~Hamilton \cite{Hamilton2} showed that positive curvature operator is preserved by the Ricci flow in all dimensions. Moreover, Hamilton proved the following convergence theorem in dimension $4$: 

\begin{theorem}[R.~Hamilton \cite{Hamilton2}]
\label{Hamilton.dim.4}
Let $(M,g_0)$ be a compact four-manifold with positive curvature operator, and let $g(t)$, $t \in [0,T)$, denote the unique maximal solution to the Ricci flow with initial metric $g_0$. Then $T < \infty$, and the rescaled metrics $\frac{1}{6(T-t)} \, g(t)$ converge to a metric of constant sectional curvature $1$ as $t \to T$.
\end{theorem} 

It was shown by H.~Chen \cite{Chen1} that two-positive curvature operator is also preserved by the Ricci flow. As a result, Chen was able to weaken the curvature assumption in Theorem \ref{Hamilton.dim.4}:

\begin{theorem}[H.~Chen \cite{Chen1}]
\label{chen.thm}
Let $(M,g_0)$ be a compact four-manifold with two-positive curvature operator, and let $g(t)$, $t \in [0,T)$, denote the unique maximal solution to the Ricci flow with initial metric $g_0$. Then $T < \infty$, and the rescaled metrics $\frac{1}{6(T-t)} \, g(t)$ converge to a metric of constant sectional curvature $1$ as $t \to T$.
\end{theorem}  

Other convergence theorems in dimension $4$ were established by B.~Andrews and H.~Nguyen \cite{Andrews-Nguyen} and by C.~Margerin \cite{Margerin2}.

The Ricci flow on manifolds of arbitrary dimension was first studied by G.~Huisken \cite{Huisken}. To describe this result, we consider a Riemannian manifold $(M,g_0)$ of dimension $n \geq 4$. We may decompose the curvature tensor of $(M,g_0)$ as 
\[R_{ijkl} = U_{ijkl} + V_{ijkl} + W_{ijkl},\] 
where $U_{ijkl}$ denotes the part of the curvature tensor associated with the scalar curvature, $V_{ijkl}$ is the part of the curvature tensor associated with the tracefree Ricci curvature, and $W_{ijkl}$ denotes the Weyl tensor. 

\begin{theorem}[G.~Huisken \cite{Huisken}]
\label{huisken.thm}
Let $(M,g_0)$ be a compact Riemannian manifold of dimension $n \geq 4$ with positive scalar curvature. Suppose that the curvature tensor of $(M,g_0)$ satisfies the pinching condition 
\[|V|^2 + |W|^2 < \delta(n) \, |U|^2,\] 
where $\delta(4) = \frac{1}{5}$, $\delta(5) = \frac{1}{10}$, and 
\[\delta(n) = \frac{2}{(n-2)(n+1)}\] 
for $n \geq 6$. Finally, let $g(t)$, $t \in [0,T)$, denote the unique maximal solution to the Ricci flow with initial metric $g_0$. Then $T < \infty$, and the rescaled metrics $\frac{1}{2(n-1)(T-t)} \, g(t)$ converge to a metric of constant sectional curvature $1$ as $t \to T$. 
\end{theorem}

In particular, Theorem \ref{huisken.thm} implies that the Ricci flow evolves any initial metric with sufficiently pinched curvature to a constant curvature metric. Similar results in this direction were obtained by C.~Margerin \cite{Margerin1} and S.~Nishikawa \cite{Nishikawa}. 

Theorem \ref{chen.thm} was generalized to higher dimensions in \cite{Bohm-Wilking}: 

\begin{theorem}[C.~B\"ohm, B.~Wilking \cite{Bohm-Wilking}]
\label{bw}
Let $(M,g_0)$ be a compact Riemannian manifold of dimension $n \geq 4$ with two-positive curvature operator, and let $g(t)$, $t \in [0,T)$, denote the unique maximal solution to the Ricci flow with initial metric $g_0$. Then $T < \infty$, and the rescaled metrics $\frac{1}{2(n-1)(T-t)} \, g(t)$ converge to a metric of constant sectional curvature $1$ as $t \to T$.
\end{theorem}

The condition that the initial manifold $(M,g_0)$ has two-positive curvature operator is very restrictive in high dimensions. In the remainder of this section, we will describe the convergence theory for the Ricci flow developed in \cite{Brendle5} and \cite{Brendle-Schoen1}. This approach requires much weaker assumptions on the initial metric $g_0$. We will only sketch the main ideas; a detailed exposition can be found in the monograph \cite{Brendle-book}. 

The approach developed in \cite{Brendle5} and \cite{Brendle-Schoen1} relies in a crucial way on the notion of \textit{positive isotropic curvature}. This notion was introduced in a fundamental paper of M.~Micallef and J.D.~Moore \cite{Micallef-Moore} on the Morse index of minimal two-spheres.

\begin{definition}
Let $M$ be Riemannian manifold of dimension $n \geq 4$. We say that $M$ has nonnegative isotropic curvature if 
\[R_{1313} + R_{1414} + R_{2323} + R_{2424} - 2 \, R_{1234} \geq 0\] 
for each point $p \in M$ and every orthonormal four-frame $\{e_1,e_2,e_3,e_4\} \subset T_p M$. If the strict inequality holds, we say that $M$ has positive isotropic curvature.
\end{definition}

Micallef and Moore \cite{Micallef-Moore} proved that every compact, simply connected manifold with positive isotropic curvature is homeomorphic to $S^n$. The topology of non-simply connected manifolds with positive isotropic curvature is not completely understood. However, there are restrictions on the fundamental groups of such manifolds (cf. \cite{Fraser}).

It turns out that the Ricci flow preserves nonnegative isotropic curvature in all dimensions. This fact plays a key role in the convergence theory in higher dimensions. 

\begin{theorem}[S.~Brendle, R.~Schoen \cite{Brendle-Schoen1}]
\label{pic.is.preserved}
Let $M$ be a compact manifold of dimension $n \geq 4$, and let $g(t)$, $t \in [0,T)$, be a solution to the Ricci flow on $M$. If $(M,g(0))$ has nonnegative isotropic curvature, then $(M,g(t))$ has nonnegative isotropic curvature for all $t \in [0,T)$.
\end{theorem}

The proof of Theorem \ref{pic.is.preserved} uses the maximum principle. The strategy is to show that the minimum isotropic curvature increases under the Ricci flow. This is a very delicate calculation which exploits special identities and inequalities arising from the first and second variations of the isotropic curvature. 

We next discuss some curvature conditions which are closely related to positive isotropic curvature and which are also preserved by the Ricci flow: 

\begin{theorem}
\label{product.conditions}
Let $M$ be a compact manifold of dimension $n \geq 4$, and let $g(t)$, $t \in [0,T)$, be a solution to the Ricci flow on $M$. Moreover, let $S^2(1)$ denote the two-sphere equipped with its standard metric of constant curvature $1$.
\begin{itemize}
\item[(i)] If $(M,g(0)) \times \mathbb{R}$ has nonnegative isotropic curvature, then $(M,g(t)) \times \mathbb{R}$ has nonnegative isotropic curvature for all $t \in [0,T)$.
\item[(ii)] If $(M,g(0)) \times \mathbb{R}^2$ has nonnegative isotropic curvature, then $(M,g(t)) \times \mathbb{R}^2$ has nonnegative isotropic curvature for all $t \in [0,T)$.
\item[(iii)] If $(M,g(0)) \times S^2(1)$ has nonnegative isotropic curvature, then the product $(M,g(t)) \times S^2(1)$ has nonnegative isotropic curvature for all $t \in [0,T)$. 
\end{itemize}
\end{theorem}

The statements (i) and (ii) follow directly from Theorem \ref{pic.is.preserved}. The proof of (iii) is more subtle, and can be found in \cite{Brendle5}.

Note that the product $M \times \mathbb{R}$ has nonnegative isotropic curvature if and only if 
\[R_{1313} + \lambda^2 \, R_{1414} + R_{2323} + \lambda^2 \, R_{2424} - 2\lambda \, R_{1234} \geq 0\] 
for all orthonormal four-frames $\{e_1,e_2,e_3,e_4\}$ and all $\lambda \in [0,1]$. Moreover, the product $M \times \mathbb{R}^2$ has nonnegative isotropic curvature if and only if 
\[R_{1313} + \lambda^2 \, R_{1414} + \mu^2 \, R_{2323} + \lambda^2\mu^2 \, R_{2424} - 2\lambda\mu \, R_{1234} \geq 0\] 
for all orthonormal four-frames $\{e_1,e_2,e_3,e_4\}$ and all $\lambda,\mu \in [0,1]$. Finally, $M \times S^2(1)$ has nonnegative isotropic curvature if and only if 
\begin{align*} 
&R_{1313} + \lambda^2 \, R_{1414} + \mu^2 \, R_{2323} + \lambda^2\mu^2 \, R_{2424} - 2\lambda\mu \, R_{1234} \\ 
&+ (1-\lambda^2)(1-\mu^2) \geq 0 
\end{align*}
for all orthonormal four-frames $\{e_1,e_2,e_3,e_4\}$ and all $\lambda,\mu \in [0,1]$. 

Using Theorem \ref{product.conditions}, one can prove the following convergence theorem for the Ricci flow:

\begin{theorem}[S.~Brendle \cite{Brendle5}]
\label{convergence.ricci.flow}
Let $(M,g_0)$ be a compact Riemannian manifold of dimension $n \geq 4$ with the property that 
\[R_{1313} + \lambda^2 \, R_{1414} + R_{2323} + \lambda^2 \, R_{2424} - 2\lambda \, R_{1234} > 0\]
for all orthonormal four-frames $\{e_1,e_2,e_3,e_4\}$ and all $\lambda \in [0,1]$. Let $g(t)$, $t \in [0,T)$, denote the unique maximal solution to the Ricci flow with initial metric $g_0$. Then $T < \infty$, and the rescaled metrics $\frac{1}{2(n-1)(T-t)} \, g(t)$ converge to a metric of constant sectional curvature $1$ as $t \to T$. 
\end{theorem}

Theorem \ref{convergence.ricci.flow} contains many known convergence results as subcases. In particular, Theorem \ref{huisken.thm} and Theorem \ref{bw} are both special cases of Theorem \ref{convergence.ricci.flow}. Moreover, Theorem \ref{convergence.ricci.flow} applies to any initial manifold $(M,g_0)$ with $1/4$-pinched sectional curvatures.

\begin{corollary}[S.~Brendle, R.~Schoen \cite{Brendle-Schoen1}] 
\label{pinching}
Let $(M,g_0)$ be a compact Riemannian manifold of dimension $n \geq 4$. We assume that $(M,g_0)$ is strictly $1/4$-pinched in the pointwise sense, so that $0 < K(\pi_1) < 4 \, K(\pi_2)$ for all points $p \in M$ and all two-planes $\pi_1,\pi_2 \subset T_p M$. Let $g(t)$, $t \in [0,T)$, denote the unique maximal solution to the Ricci flow with initial metric $g_0$. Then $T < \infty$, and the rescaled metrics $\frac{1}{2(n-1)(T-t)} \, g(t)$ converge to a metric of constant sectional curvature $1$ as $t \to T$. 
\end{corollary}

In particular, Corollary \ref{pinching} implies that any Riemannian manifold with strictly $1/4$-pinched curvature is diffeomorphic to a spherical space form. This result is known as the Differentiable Sphere Theorem: 

\begin{diffeo.sph.thm}[S.~Brendle, R.~Schoen \cite{Brendle-Schoen1}]
Suppose that $M$ is a compact Riemannian manifold which is strictly $1/4$-pinched in the pointwise sense, so that $0 < K(\pi_1) < 4 \, K(\pi_2)$ for all points $p \in M$ and all two-planes $\pi_1,\pi_2 \subset T_p M$. Then $M$ is diffeomorphic to a quotient $S^n / \Gamma$, where $\Gamma$ is a finite group of isometries acting freely.
\end{diffeo.sph.thm}

The Sphere Theorem has a long history, dating back to a question posed by H.~Hopf in the 1940s (see \cite{Brendle-Schoen3} for a survey). The classical Sphere Theorem of Berger \cite{Berger} and Klingenberg \cite{Klingenberg} asserts that any compact, simply connected Riemann manifold whose sectional curvatures all lie in the interval $(1,4]$ is homeomorphic to the sphere. The Differentiable Sphere Theorem improves this result in several respects. Most importantly, it provides a classification up to diffeomorphism (not just up to homeomorphism). Moreover, the Differentiable Sphere Theorem applies to non-simply connected manifolds, and only requires a pointwise pinching assumption. 

Using the strict maximum principle, one can analyze the borderline case in Theorem \ref{convergence.ricci.flow} (see \cite{Brendle-book}, Theorem 9.33).

\begin{theorem}[S.~Brendle \cite{Brendle-book}]
\label{borderline}
Let $(M,g_0)$ be a compact, simply connected Riemannian manifold of dimension $n \geq 4$ with the property that 
\[R_{1313} + \lambda^2 \, R_{1414} + R_{2323} + \lambda^2 \, R_{2424} - 2\lambda \, R_{1234} \geq 0\] 
for all orthonormal four-frames $\{e_1,e_2,e_3,e_4\}$ and all $\lambda \in [0,1]$. Moreover, let $g(t)$, $t \in [0,T)$, denote the unique maximal solution to the Ricci flow with initial metric $g_0$. Then one of the following statements holds: 
\begin{itemize}
\item[(i)] $T < \infty$, and the rescaled metrics $\frac{1}{2(n-1)(T-t)} \, g(t)$ converge to a metric of constant sectional curvature $1$ as $t \to T$.
\item[(ii)] $(M,g_0)$ is a K\"ahler manifold.
\item[(iii)] $(M,g_0)$ is isometric to a symmetric space.
\item[(iv)] $(M,g_0)$ is isometric to a product of two manifolds of lower dimension.
\end{itemize}
\end{theorem}

If $(M,g_0)$ is a Riemannian manifold with weakly $1/4$-pinched sectional curvatures, then the curvature assumption in Theorem \ref{borderline} is satisfied. Hence, we can draw the following conclusion:

\begin{corollary}[S.~Brendle, R.~Schoen \cite{Brendle-Schoen2}] 
Let $(M,g_0)$ be a compact, simply connected manifold of dimension $n \geq 4$. Moreover, suppose that $(M,g_0)$ is weakly $1/4$-pinched in the pointwise sense, so that $0 \leq K(\pi_1) \leq 4 \, K(\pi_2)$ for all points $p \in M$ and all two-planes $\pi_1,\pi_2 \subset T_p M$. Let $g(t)$, $t \in [0,T)$, be the unique maximal solution to the Ricci flow with initial metric $g_0$. Then one of the following statements holds: 
\begin{itemize}
\item[(i)] $T < \infty$, and the rescaled metrics $\frac{1}{2(n-1)(T-t)} \, g(t)$ converge to a metric of constant sectional curvature $1$ as $t \to T$.
\item[(ii)] $(M,g_0)$ is isometric to a symmetric space. 
\end{itemize}
\end{corollary}

There is a well-known classification of all symmetric spaces with positive sectional curvature. The list of such spaces includes the following examples: 
\begin{itemize} 
\item The round sphere $S^n$.
\item The complex projective space $\mathbb{CP}^m$. 
\item The quaternionic projective space $\mathbb{HP}^m$. 
\item The projective plane over the octonions.
\end{itemize}
Each of these spaces has weakly $1/4$-pinched sectional curvatures.

Finally, we note that the ideas described above can be used to improve many classical results on the Ricci flow. For example, R.~Hamilton \cite{Hamilton4} showed that any solution to the Ricci flow with nonnegative curvature operator satisfies the Harnack inequality 
\[\frac{\partial}{\partial t} R + \frac{1}{t} \, R + 2 \, \partial_i R \, v^i + 2 \, \text{\rm Ric}_{ij} \, v^i \, v^j \geq 0.\] 
Here, $R$ denotes the scalar curvature, $\text{\rm Ric}$ denotes the Ricci tensor of the evolving metric, and $v$ is an arbitrary tangent vector. In \cite{Brendle6}, we showed that Hamilton's Harnack inequality holds under a much weaker curvature assumption; in fact, it is enough to assume that $M \times \mathbb{R}^2$ has nonnegative isotropic curvature. The techniques discussed above also yield new results about ancient solutions to the Ricci flow (see \cite{Brendle-Huisken-Sinestrari} for details).

\end{document}